\documentclass[11pt]{article}
\pdfoutput=1
\usepackage{amsbsy,amsmath,amsthm,amssymb,graphicx,verbatim}
\usepackage{setspace, float, subfig}
\usepackage[longnamesfirst,sort]{natbib}

\usepackage{multirow}

\usepackage{geometry}
\newtheorem{theorem}{Theorem}

\newtheorem{proposition}{Proposition}
\theoremstyle{remark}
\newtheorem{remark}{Remark}

\newcommand{\sX}{\mathsf{X}}
\newfont{\msbm}{msbm10 at 11pt}
\newcommand {\R} {\mathbb{R}}

\newcommand {\N} {\mbox{\msbm N}}

\newcommand{\Var}[0]{\text{Var}}

\begin{document}
\onehalfspacing

\title{Relative fixed-width stopping rules for Markov chain Monte Carlo simulations}

\author{James M. Flegal \\ Department of Statistics \\ University of California, Riverside \\ {\tt jflegal@ucr.edu} \and Lei Gong \\ Department of Statistics \\ University of California, Riverside \\ {\tt lei.gong@email.ucr.edu} }  

\date{\today}

\maketitle

\begin{abstract}
Markov chain Monte Carlo (MCMC) simulations are commonly employed for estimating features of a target distribution, particularly for Bayesian inference.  A fundamental challenge is determining when these simulations should stop.  We consider a sequential stopping rule that terminates the simulation when the width of a confidence interval is sufficiently small relative to the size of the target parameter.  Specifically, we propose relative magnitude and relative standard deviation stopping rules in the context of MCMC.  In each setting, we develop sufficient conditions for asymptotic validity, that is conditions to ensure the simulation will terminate with probability one and the resulting confidence intervals will have the proper coverage probability.  Our results are applicable in a wide variety of MCMC estimation settings, such as expectation, quantile, or simultaneous multivariate estimation.  Finally, we investigate the finite sample properties through a variety of examples and provide some recommendations to practitioners.  

\end{abstract}

\section{Introduction} \label{sec:intro}

Markov chain Monte Carlo (MCMC) methods allow exploration of intractable probability distributions by constructing a Markov chain whose stationary distribution equals the desired distribution.  A major challenge for practitioners is determining how long to run an MCMC simulation.  Many experiments employ a fixed-time rule to terminate the simulation; that is, the procedure terminates after $n$ iterations, where $n$ is determined heuristically.  Indeed, some simulations are so complex that this is the only practical approach, but that is not so for most experiments.  

Alternatively, many practitioners use convergence diagnostics to determine if $n$ is sufficiently large \citep[for a review see][]{cowl:carl:1996}.  Although practical, these methods are mute about the quality of the resulting estimates \citep{fleg:hara:jone:2008}.  Moreover, they can introduce bias directly in to the estimates \citep{cowl:robe:rose:1999}.  

We instead advocate terminating the simulation when an estimate is sufficiently accurate for the analytic purpose that motivates the inquiry.  In other words, the simulation is terminated the first time a confidence interval width for a desired quantity is sufficiently small.  We refer to such a procedure as a sequential fixed-width stopping rule and note the total simulation effort will be random.  

As we show later, fixed-width methods are especially desirable because they are theoretically justified and constrained by few assumptions.  The simplest fixed-width rule, first studied in MCMC by \cite{jone:hara:caff:neat:2006}, stops the simulation when the width of a confidence interval based on an ergodic average is less than a user-specified value, say $\epsilon$.  \cite{fleg:hara:jone:2008} and \cite{jone:hara:caff:neat:2006} show this stopping rule is superior to using convergence diagnostics as a stopping criteria.  

In this paper, we introduce relative fixed-width stopping rules that eliminate the need to specify an absolute value for $\epsilon$.  Specifically, the simulation is terminated the first time the width of a confidence interval is sufficiently small relative to the \textit{size} of a target parameter.  We consider two measures of size, magnitude and standard deviation.  Further, we illustrate the utility of relative fixed-width stopping rules for simultaneous estimation of multiple parameters. 




Specificity requires some notation.  Let $\pi$ denote a probability distribution having support $\sX \subseteq \mathbb{R}^{d}$, $d \ge 1$, about which we wish to make inference.  This inference is usually based on some feature of $\pi$ denoted $\theta_{\pi}$.  For example, we may want a quantile of $\pi$ or if $g: \sX \to \R$, we may need to calculate 
\[
E_{\pi} [ g(X) ] = \int_{\sX} g(x) \pi(dx) \; .
\]

Frequently $\pi$ is such that MCMC is the only viable technique for estimating $\theta_{\pi}$.  The basic MCMC method entails constructing a time-homogeneous Harris ergodic Markov chain $X = \left\{  X^{(0)} ,  X^{(1)} , \ldots\right\}$ on state space $\sX$ with $\sigma$-algebra $\mathcal{B} = \mathcal{B}(\mathsf{X})$ and invariant distribution $\pi$.  The popularity of MCMC methods result from the ease with which $X$ can be simulated \citep{robe:case:2004}.

Suppose we simulate $X$ for $n$ iterations, where $n$ is finite.  Define $\hat{\theta}_n$ as an estimator of $\theta_{\pi}$ from the observed chain.  Outside of toy examples, no matter how long our simulation, there will be an unknown Monte Carlo error, $\hat{\theta}_{n} - \theta_{\pi}$.  While it is impossible to assess this error directly, we can obtain its approximate sampling distribution if a Markov chain central limit theorem (CLT) holds.  That is, if
\begin{equation}
\label{eq:clt}
\sqrt{n} \left( \hat{\theta}_{n} - \theta_{\pi} \right) \stackrel{d}{\to} \; \text{N}(0, \sigma^{2}_{\theta})
\end{equation}
as $n \to \infty$ where $\sigma^{2}_{\theta} \in (0, \infty)$.  Denote $\lambda^2_{\theta}$ as the posterior variance associated with $\theta_{\pi}$.  Then it is important to note that due to the correlation present in a Markov chain  $\sigma^{2}_{\theta} \neq \lambda^2_{\theta}$, except in trivial cases.  

For now, suppose we have an estimator such that $\hat{\sigma}^{2}_{n} \to \sigma^{2}_{\theta}$ almost surely as $n \to \infty$.  This allows construction of a $(1-\delta)100\%$ confidence interval for $\theta_{\pi}$ with width
\begin{equation}
\label{eq:hw}
w_{\delta} = 2 z_{\delta / 2} \frac{\hat{\sigma}_{n}}{\sqrt{n}} 
\end{equation}
where $z_{\delta / 2}$ is a critical value from a standard Normal distribution.  The width at \eqref{eq:hw} allows analysts to report the uncertainty in their estimates and users to assess the practical reliability.  Moreover, we will use $w_{\delta}$ to construct sequential fixed-width stopping rules.

Our work advocates stopping the simulation the first time $w_{\delta}$ is sufficiently small.  We consider three distinct stopping rules: (i) an absolute precision rule that terminates when $w_{\delta} < \epsilon$, (ii) a relative magnitude rule that terminates when $w_{\delta} < \epsilon \left| \theta_{\pi} \right| $ and (iii) a relative standard deviation rule that terminates when $w_{\delta} < \epsilon \lambda_{\theta}$.

The theoretical properties of (i) and (ii) have been studied by \cite{glyn:whit:1992}, which we extend to establish conditions for asymptotic validity of (iii).  Asymptotic validity is important since it implies the simulation will terminate w.p.1 and the resulting confidence intervals will have the right coverage probability.  

\cite{fleg:hara:jone:2008}, \cite{fleg:jone:2010} and \cite{jone:hara:caff:neat:2006} have previously investigated (i) for MCMC expectation estimation.  We are not aware of any prior use of fixed-width methods for quantile estimation or any use of (ii) or (iii) as a stopping rule in MCMC.  The rule (iii) has significant promise in Bayesian applications since the simulation terminates the first time the length of a confidence interval is less than an $\epsilon$th fraction of the magnitude of the standard deviation of $\theta_{\pi}$.  In other words, the simulation stops when an estimate of $\theta_{\pi}$ is sufficiently accurate relative to an associated posterior standard deviation.  Another substantial benefit of rule (iii) is it easy to implement in multivariate settings since $\epsilon$ can remain constant.

There are two main assumptions for asymptotic validity.  First, we require a limiting distribution for the Monte Carlo error such as at \eqref{eq:clt}.  Second, we require a strongly consistent estimator of the associated asymptotic variance, that is $\hat{\sigma}^{2}_{n} \to \sigma^{2}_{\theta}$ almost surely as $n \to \infty$.  We later discuss these assumptions in detail for estimating expectations and quantiles.  

Finally, we investigate the finite sample properties of relative fixed-width stopping rules through three examples.  Our first example considers an independence Metropolis sampler to explore an exponential random variable.  Our second example considers exploring a mixture of bivariate Normal distributions with Metropolis Hastings and Gibbs samplers.  While these are only toy examples, we will use true parameter values to illustrate the utility of our stopping rules.  Our final example considers a Bayesian version of a logistic regression to model the presence or absence of the freshwater eel {\it Anguilla australis}.

Using these examples, we terminate the simulation with the three distinct fixed-width stopping rules and calculate confidence intervals for a vector of target parameters.  Over replicated simulations, all the finite sample empirical coverage probabilities are close to a specified nominal level.  Thus, fixed-width stopping rules provide a theoretically valid and practically accurate procedure to determine when to stop a MCMC simulation.

For Bayesian practitioners, we advocate the relative standard deviation fixed-width stopping rule (iii) since it is easy to implement and applicable in multivariate settings without a priori knowledge of the target parameter size.  As our examples show, setting $\epsilon = 0.02$ provides excellent results in a wide variety of univariate and multivariate settings.  

The rest of this paper is organized as follows.  Section~\ref{sec:fixed} formally introduces relative fixed-width stopping rules and establishes asymptotic validity.  Section~\ref{sec:exp} investigates fixed-width stopping procedures when estimating expectations and quantiles.  Section~\ref{sec:examples} studies the finite sample properties in three numerical examples and concludes with a discussion that provides some recommendations to practitioners.

\section{Sequential fixed-width procedures} \label{sec:fixed}
In this section, we obtain conditions that ensure asymptotic validity of fixed-width procedures.  The primary assumption is a limiting distribution holds for the Monte Carlo error such as at \eqref{eq:clt}.  Actually, we require the limiting process at \eqref{eq:clt} satisfy a stronger condition, in particular, a functional central limit theorem (FCLT).  

Define an interval 
\[
C[n] = \left( \hat{\theta}_{n} - z_{\delta / 2} \hat{\sigma}_{n} / \sqrt{n} \text{ , } \hat{\theta}_{n} + z_{\delta / 2} \hat{\sigma}_{n} / \sqrt{n} \right) \; .
\]
If \eqref{eq:clt} holds and $\hat{\sigma}_{n}$ is weakly consistent for $\sigma_{\theta}$, then $C[n]$ achieves the nominal coverage level as the sample size $n \to \infty$.  Thus we have a valid confidence interval provided the sample size is permitted to go to $\infty$.  

Now consider a sequential procedure that terminates the simulation when the length of a confidence interval drops below a prescribed level $\epsilon$.  We will refer to this type of stopping rule as an absolute precision fixed-width stopping rule.  For such a rule, the time at which the simulation terminates is defined by
\[
\tilde{T} (\epsilon) = \inf \left\{ n \ge 0 : 2 z_{\delta / 2} \hat{\sigma}_{n} / \sqrt{n} \le \epsilon \right\} \; .
\]
Unfortunately, use of this stopping rule is insufficient because $\tilde{T} (\epsilon)$ can terminate much too early if $\hat{\sigma}_{n}$ is poorly behaved for small $n$ \citep{glyn:whit:1992}.  Instead, suppose $p(n)$ is a positive function that decreases monotonically such that $p(n) = o(n^{-1/2})$ as $n \to \infty$ and let $n^{*}$ be the desired minimum simulation effort (a reasonable default is $p(n) = \epsilon I(n \le n^{*}) + n^{-1} $).  Then an absolute precision stopping rule terminates the simulation at
\[
T_1 (\epsilon) = \inf \left\{ n \ge 0 : 2 z_{\delta / 2} \hat{\sigma}_{n} / \sqrt{n} + p(n) \le \epsilon \right\} \; .
\]

The following result, an immediate consequence of Theorem~1 in \cite{glyn:whit:1992}, yields asymptotic validity of the sequential stopping rule $T_1 (\epsilon)$.  Note the desired coverage probability will be obtained in an asymptotic sense as $\epsilon \to 0$.
\begin{proposition} \label{thm:gw1}
Suppose a FCLT at \eqref{eq:clt} holds.  If $\hat{\sigma}_{n} \to \sigma_{\theta}$ w.p.1 as $n \to \infty$, then as $n \to \infty$ or $\epsilon \to 0$ the simulation will terminate w.p.1 and 
\begin{equation*}
Pr\left( \theta_{\pi} \in C[T_1 (\epsilon)] \right) \to 1 - \delta \; .
\end{equation*}
\end{proposition}

\begin{remark}
\cite{glyn:whit:1992} show weak consistency of $\hat{\sigma}_{n}$ is not enough to ensure asymptotic validity .
\end{remark}

The stopping rule $T_1 (\epsilon)$ has previously been used for estimating expectations in MCMC \citep{fleg:jone:2010, fleg:hara:jone:2008, jone:hara:caff:neat:2006}.  We further show this rule works well for MCMC estimation of quantiles in the following section.  The challenge in both settings is finding a strongly consistent estimator of $\sigma_{\theta}$.

One can consider a variant of the stopping rule $T_1 (\epsilon)$ known as a relative precision stopping rule, which avoids having to choose an absolute value for $\epsilon$.  Simply put, the simulation is run until the length of a confidence interval is less than an $\epsilon$th fraction of the magnitude of the parameter of interest, $\theta_{\pi}$.  Using $\hat{\theta}_{n}$ as an estimator of $\theta_{\pi}$ yields the following relative magnitude stopping rule
\[
T_2 (\epsilon) = \inf \left\{ n \ge 0 : 2 z_{\delta / 2} \hat{\sigma}_{n} / \sqrt{n} + p(n) \le \epsilon \left| \hat{\theta}_{n} \right| \right\} \; .
\]
For large $n$, $T_2 (\epsilon)$ will behave like $T_1 (\epsilon |\theta_{\pi}| )$.  The following obtains asymptotic validity of $T_2 (\epsilon)$, which is a direct consequence of Theorem~3 in \cite{glyn:whit:1992}.
\begin{proposition} \label{thm:gw3}
Suppose a FCLT at \eqref{eq:clt} holds and $|\theta_{\pi}| > 0$.  If $\hat{\theta}_{n} \to \theta_{\pi}$ w.p.1 and $\hat{\sigma}_{n} \to \sigma_{\theta}$ w.p.1 as $n \to \infty$, then as $n \to \infty$ or $\epsilon \to 0$ the simulation will terminate w.p.1 and 
\begin{equation*}
Pr\left( \theta_{\pi} \in C[T_2 (\epsilon)] \right) \to 1 - \delta \; .
\end{equation*}
\end{proposition}

Note that Proposition~\ref{thm:gw3} requires $\hat{\theta}_{n} \to \theta_{\pi}$ w.p.1 along with necessary conditions of Proposition~\ref{thm:gw1}.  In general stochastic simulations, this condition does not immediately follow from \eqref{eq:clt} \citep[see Example 2 of][]{glyn:whit:1988} but is readily available when $\theta_{\pi}$ is an expectation via the Markov chain strong law of large numbers (SLLN).  

While $T_2 (\epsilon)$ has some support in the operations research literature, it makes little intuitive sense in Bayesian settings.  Specifically, if $\theta_{\pi} = 0$ then $T_2 (\epsilon)$ will be theoretically invalid and poorly behaved in finite simulations.  In addition, $T_2 (\epsilon)$ could be problematic even when $\theta_{\pi}$ is close to zero, which we illustrate through example in Section~\ref{sec:examples}.

Given the popularity of MCMC in Bayesian settings, it is useful to consider another specifically designed variant of $T_1 (\epsilon)$.  To this end, we propose a stopping rule that terminates the simulation when the length of a confidence interval is less than an $\epsilon$th fraction of the magnitude of $\lambda_{\theta}$, i.e.\ the posterior standard deviation of $\theta_{\pi}$.  Suppose $\hat{\lambda}_{n}$ is an estimator of $\lambda_{\theta}$ and consider the following stopping rule
\[
T_3 (\epsilon) = \inf \left\{ n \ge 0 : 2 z_{\delta / 2} \hat{\sigma}_{n} / \sqrt{n} + p(n) \le \epsilon \hat{\lambda}_{n} \right\} \; .
\]
For large $n$, $T_3 (\epsilon)$ will behave like $T_1 (\epsilon \lambda_{\theta} )$.  The benefit of using $T_3 (\epsilon)$ is that $\epsilon$ is selected as a fraction rather than in the units of the target parameter.  Hence, a single value of $\epsilon$ would be appropriate for target parameters of any magnitude.  Naturally, decreasing $\epsilon$ would decrease the uncertainty of the resulting estimates and could be done simultaneously for multiple parameters.  The following establishes asymptotic validity of $T_3 (\epsilon)$, which we prove in Appendix~\ref{app:T3}.
\begin{theorem} \label{thm:T3}
Suppose a FCLT at \eqref{eq:clt} holds and $\lambda_{\theta}>0$.  If $\hat{\lambda}_{n} \to \lambda_{\theta}$ w.p.1 and $\hat{\sigma}_{n} \to \sigma_{\theta}$ w.p.1 as $n \to \infty$, then as $n \to \infty$ or $\epsilon \to 0$ the simulation will terminate w.p.1 and 
\begin{equation*}
Pr\left( \theta_{\pi} \in C[T_3 (\epsilon)] \right) \to 1 - \delta \; .
\end{equation*}
\end{theorem}

Note the only additional condition required for Theorem~\ref{thm:T3} is a strongly consistent estimator of $\lambda_{\theta}$.  For expectations, an estimator is readily available via the Markov chain SLLN.  In the case of quantiles, we discuss a viable estimator in the following section.

The benefit of the stopping rule $T_3 (\epsilon)$ is twofold.  First, one only needs to specify a relative $\epsilon$, and hence no knowledge about the magnitude is required.  Second, when estimating multiple parameters a single $\epsilon$ will suffice to obtain estimates whose uncertainty will be comparable relative to their standard deviations.  In other words, we have developed a simple, yet informative, stopping criteria applicable in multivariate settings.  In these settings, one could address the issue of multiplicity by adjusting the critical value appropriately.  We illustrate this procedure via examples in Section~\ref{sec:examples}, and show the resulting simultaneous confidence regions obtain at least the nominal coverage probability.

\begin{remark}
Asymptotic validity of relative stopping rules can be established when \eqref{eq:clt} is replaced by a more general $\R$-valued stochastic process \citep{glyn:whit:1992}.  The generalization enables consideration of $\theta_{\pi}$ that follow non-Normal asymptotic distributions.
\end{remark}


\section{Applications} \label{sec:exp}
This section demonstrates that fixed-width stopping rules are appropriate for MCMC estimation of expectations and quantiles.  This is an important contribution since we know of no other formal stopping criteria applicable in both settings.  \cite{raft:lewi:1992} propose a heuristic approach to terminating an MCMC simulation when the primary interest is quantile estimation.  However, \cite{broo:robe:1999} argue ``in the case where quantiles themselves are not of interest, this method should be used with caution''.  

First, we require a bit more notation to describe sufficient mixing conditions for a Markov chain CLT and consistent estimation of the asymptotic variance.  An interested reader is directed to \cite{meyn:twee:glyn:2009} and \cite{robe:rose:2004} for more on Markov chain theory.  

Recall $X$ is a Harris ergodic Markov chain on state space $\sX$ with $\sigma$-algebra $\mathcal{B} = \mathcal{B}(\mathsf{X})$ and invariant distribution $\pi$.  Denote the $n$-step Markov kernel associated with $X$ as $P^{n}(x,dy)$ for $n \in \N$.  Then if $A \in \mathcal{B}(\sX)$ and $k\in \{0, 1, 2, \ldots \}$, $P^{n}(x,A) = \Pr (X_{k+n} \in A | X_{k} =x) $.  Let $\|\cdot\|$ denote the total variation norm.  Let $M: \sX \mapsto \R^+$ and $\gamma: \N \mapsto \R^+$ be decreasing such that
\begin{equation} \label{eq:tvn}
\Vert P^n (x, \cdot) - \pi(\cdot) \Vert \leq M(x) \gamma(n) \; .
\end{equation}
\textit{Polynomial ergodicity of order $m$} where $m \ge 0$ means \eqref{eq:tvn} holds with $E_{\pi} M < \infty$ and $\gamma(n) = n^{-m}$ for all $X_{0} = x$.  \textit{Geometrical ergodicity} means \eqref{eq:tvn} holds with $\gamma(n) = t^n$ for some $0 < t < 1$ for all $X_{0} = x$.  \textit{Uniform ergodicity} means \eqref{eq:tvn} holds with $M$ bounded and $\gamma(n) = t^n$ for some $0 < t < 1$.

Establishing \eqref{eq:tvn} directly can be challenging, but some constructive techniques are available \citep{jarn:robe:2002, meyn:twee:glyn:2009}.  Most literature on MCMC algorithms focuses on establishing geometric and uniform ergodicity, see e.g.\ \citet{hobe:2011}, \citet{jone:hobe:2001}, \citet{john:jone:neat:2011}, \citet{meng:twee:1996}, \citet{robe:twee:1996} and \citet{tier:1994}.  Less has been said concerning polynomial ergodicity, but an interested reader is directed to \citet{douc:fort:moul:soul:2004}, \citet{fort:moul:2000}, \citet{fort:moul:2003}, \citet{jarn:robe:2002}, \citet{jarn:robe:2007} and \citet{jarn:twee:2003}.


\subsection{Expectations}
For MCMC estimation of an expectation, one can obtain all the necessary conditions for asymptotic validity of fixed-width stopping rules.  Let $g : \sX \to \R$, then we consider estimation of
\[
\mu_{g} := E_{\pi} [ g(X) ] = \int_{\sX} g(x) \pi(dx) \; .
\]
Estimating $\mu_{g}$ is natural by appealing a Markov chain SLLN, a special case of the Birkhoff Ergodic Theorem \citep[p.\ 558][]{fris:gray:1997}.  Specifically, if $ E_{\pi} |g| < \infty$ then w.p.1
\begin{equation*}
\bar{g}_{n} := \frac{1}{n} \sum_{i=0}^{n-1} g ( X^{(i)} ) \rightarrow \mu_{g} \text{ as } n \rightarrow \infty \; .
\end{equation*}
Hence the SLLN yields strongly consistent estimators of $\mu_{g}$ and $\lambda^2_{\theta} = \Var [ g ]$ (provided $ E_{\pi} g^2 < \infty$) necessary for Proposition~\ref{thm:gw3} and Theorem~\ref{thm:T3}, respectively.  

We can obtain an approximate sampling distribution for the Monte Carlo error via a Markov chain CLT if
\begin{equation}
\label{eq:mu clt}
\sqrt{n} (\bar{g}_{n} - \mu_{g}) \stackrel{d}{\to} \; \text{N}(0, \sigma^{2}_{g})
\end{equation}
as $n \to \infty$ where $\sigma^{2}_{g} \in (0, \infty)$.  Conditions that ensure \eqref{eq:mu clt} can be found in \citet{chan:geye:1994}, \citet{jone:2004}, \citet{meyn:twee:glyn:2009}, \citet{robe:rose:2004}  and \citet{tier:1994}.  For example, if $X$ is geometrically ergodic and $E_{\pi} |g|^{2+\epsilon} < \infty$ for some $\epsilon >0$, then \eqref{eq:mu clt} holds.  Fortunately, Markov chains frequently enjoy a FCLT under the same conditions \citep{ooda:yosh:1972, ibra:1962}.  

There are many strongly consistent variance estimation techniques applicable for $\sigma^{2}_{g}$ in MCMC settings including batch means \citep{fleg:jone:2010, jone:hara:caff:neat:2006}, spectral variance techniques \citep{fleg:jone:2010} and regenerative simulation \citep{hobe:jone:pres:rose:2002, mykl:tier:yu:1995}.  We consider only non-overlapping batch means (BM) because it is easy to implement and available in many software packages, e.g.\ the {\tt mcmcse} package available on {\tt CRAN}.

In BM the output is broken into $a_n$ batches where each batch is $b_n$ iterations in length.  Suppose the algorithm is run for a total of $n=a_n b_n$ iterations and define
\begin{equation*}
\bar{Y}_{j} := \frac{1}{b_n} \sum_{i=(j-1)b_n+1}^{jb_n} g (X_{i}) \hspace*{5mm} \text{ for } j=1,\ldots,a_n \; .
\end{equation*}
The BM estimate of $\sigma^{2}_{g}$ is
\begin{equation}
\label{eq:bmvar}
\hat{\sigma}_{n}^{2} = \frac{b_n}{a_n-1} \sum_{j=1}^{a_n} (\bar{Y}_{j} - \bar{g}_n)^{2}  \; .
\end{equation}

In general, the BM estimator at \eqref{eq:bmvar} is not a consistent estimator of $\sigma^{2}_{g}$.  However, \citet{jone:hara:caff:neat:2006} establish necessary conditions for $\hat{\sigma}_{n}^{2} \to \sigma^{2}_{g}$ with probability 1 as $n \to \infty$ if the batch size and the number of batches are allowed to increase as the overall length of the simulation increases.  Setting $b_{n} = \lfloor n^{\tau} \rfloor$ and $a_{n} = \lfloor n / b_{n} \rfloor$, the regularity conditions require that $X$ be geometrically ergodic, $E_{\pi} |g|^{2+\epsilon_{1}+\epsilon_{2}} < \infty$ for some $\epsilon_{1} >0$, $\epsilon_{2} >0$ and $(1+\epsilon_{1}/2)^{-1} < \tau < 1$.  A common choice of $\tau=1/2$ (i.e., $b_{n}=\lfloor \sqrt{n} \rfloor$ and $a_{n} = \lfloor n/ b_{n} \rfloor$) has been shown to work well in applications \citep{jone:hara:caff:neat:2006, fleg:jone:2010, fleg:hara:jone:2008}.  

\begin{remark}
Most sampling plans require storing the entire Markov chain to allow for recalculations as the batch size increases with $n$.  If storage is a concern, one could consider increasing the batch size of the form $b_n \in \{2, 4, 8, ... , 2^k , ...\}$ in an effort to reduce memory usage.  One can establish strong consistency for the BM variance estimator with such a sampling plan using results in \cite{jone:hara:caff:neat:2006} and \cite{bedn:latu:2007}.
\end{remark}



\subsection{Quantiles} \label{sec:other}
It is routine to estimate univariate quantiles associated with $\pi$, especially in Bayesian applications.  To this end, let $W \sim \pi$ and recall $g : \sX \to \R$.  Setting $V=g(W)$, we consider estimation of the quantiles associated with the univariate distribution of $V$.  Suppose $F_{V}$ denotes the cumulative distribution function of $V$, then our goal is to obtain 
\begin{equation*}
\xi_{q} := F_{V}^{-1} (q) = \inf \{ v: F_{V}(v) \ge q \} \; .
\end{equation*}
Little has been formally said regarding MCMC estimation of quantiles, but we outline the current state of understanding \citep[for more details see][]{fleg:jone:neat:2012}.

A natural estimator of $\xi_{q}$ is the inverse of the empirical distribution function given by
\begin{equation} 
\label{eq:estimator} 
\hat{\xi}_{n,q} : = Y_{n(j+1)} \quad \text{ where } ~~ j \le nq < j + 1\; ,
\end{equation}
where $Y_{n(j)}$ denotes the $j$th order statistic of $\{Y_{0}, \ldots, Y_{n-1}\} = \{ g(X_{0}), \ldots, g(X_{n-1}), \}$.  If $X$ is Harris recurrent and $\xi_{q}$ is the unique solution $y$ of $F_{V} (y -) \le q \le F_{V} (y)$, then $\hat{\xi}_{n,q} \to \xi_{q}$ w.p.1 as $n \rightarrow \infty$ \citep{fleg:jone:neat:2012}.

Under stronger mixing conditions on $X$, one can obtain a Markov chain CLT.  To this end, define
\begin{equation*}
\sigma^{2}(y) := \text{Var}_{\pi} \left[ I(Y_{0} \le y) \right] + 2 \sum_{k=1}^{\infty} \text{Cov}_{\pi} \left[ I(Y_{0} \le y), I(Y_{k} \le y) \right] \; . 
\end{equation*}
Suppose there is $\epsilon > 0$ such that $X$ is polynomially ergodic of order $2.5 + \epsilon$.  If $F_{V}$ has a density $f_{V}$ positive and bounded in some neighborhood of $\xi_{q}$, then as $n \to \infty$ 
\begin{equation}
\label{eq:CLT q}
\sqrt{n}(\hat{\xi}_{n,q} - \xi_{q}) \stackrel{d}{\to} \text{N}(0, \gamma^{2}(\xi_{q})) \;,
\end{equation}
where $\gamma^{2}(\xi_{q}) =  \sigma^{2} (\xi_{q}) / [f_{V}(\xi_{q})]^{2}$, provided $\sigma^{2} (\xi_{q}) > 0$ \citep{fleg:jone:neat:2012}.  A FCLT is extremely likely to hold under similar conditions, but we are unaware of a formal proof.

Estimation of the variance from the asymptotic Normal distribution at \eqref{eq:CLT q} is broken into two parts.  First, we plug in $\hat{\xi}_{n,q}$ for $\xi_{q}$ and separately consider estimating $f_{V}(\hat{\xi}_{n,q})$ and $\sigma^{2} (\hat{\xi}_{n,q})$.  Estimating $f_{V}(\hat{\xi}_{n,q})$ uses a kernel density approach with a gaussian kernel, which we denote as $\hat{f}_{V}(\hat{\xi}_{n,q})$.  There are well known conditions guaranteeing strongly consistent estimation of the density at a point \citep[see e.g.][]{kim:lee:2005, yu:1993}.

We will use BM for estimating $\sigma^{2}(\hat{\xi}_{n,q})$.  Suppose we have $n=a_n b_n$ iterations, then for $k=0,\ldots,a_n-1$ define $\bar{U}_{k}(\hat{\xi}_{n,q}) := b_n^{-1} \sum_{i=0}^{b_n-1} I(Y_{k b_n +i} \le \hat{\xi}_{n,q}) $.  The BM estimate of $\sigma^{2}(\hat{\xi}_{n,q})$ is
\begin{equation*}
\hat{\sigma}_{BM}^{2} (\hat{\xi}_{n,q}) = \frac{b_n}{a_n-1} \sum_{k=0}^{a_n-1} \left( \bar{U}_{k}(\hat{\xi}_{n,q}) - \bar{U}_n(\hat{\xi}_{n,q}) \right)^{2}  \; . 
\end{equation*}

Combining $\hat{f}_{V}(\hat{\xi}_{n,q})$ and $\hat{\sigma}_{BM}^{2} (\hat{\xi}_{n,q})$, we estimate $\gamma^{2}(\xi_{q})$ with
\[
\hat{\gamma}^{2}(\hat{\xi}_{n,q}) := \frac {\hat{\sigma}_{BM}^{2} (\hat{\xi}_{n,q})}{\hat{f}_{V}(\hat{\xi}_{n,q})} \; .
\]
This approach is implemented in the R package {\tt mcmcse} which is used to perform the computations in our examples.  \cite{fleg:jone:neat:2012} outline the conditions that ensure strong consistency of this estimator.

The relative standard deviation fixed-width stopping rule of Theorem~\ref{thm:T3} requires a estimation of 
\[
\lambda_{\theta} = \frac{q(1-q)}{ f_{V}(\xi_{q}) } \;.
\]
We use the same kernel density estimate resulting in
\[
\hat{\lambda}_{n} =  \frac{q(1-q)}{ \hat{f}_{V}(\hat{\xi}_{n,q}) } \;.
\]

\section{Numerical studies} \label{sec:examples}
This section investigates the finite sample properties of fixed-width stopping rules through a variety of simulations.   In each example, we independently repeat the MCMC simulation to evaluate the resulting finite sample confidence intervals.  Naturally, this evaluation requires the true parameter values.  In our first two examples, the true values are readily available.  In our final example, the truth was estimated using an independent long run of the MCMC sampler.  Overall, the empirical coverage probabilities obtained via fixed-width stopping rules are remarkably close to the nominal level.

Each simulation considered both expectations and quantiles with the following common methodology.  For a single replication, the same MCMC draws were used in applying the three stopping rules.  Further, we uniformly set $p(n) = \epsilon I (n < n^*) + n^{-1}$ and estimate $\sigma^2_{\theta}$ via BM methods with $b_{n}=\lfloor \sqrt{n} \rfloor$ calculated with the {\tt mcmcse} package.  Finally, standard errors for the empirical coverage probabilities equal $\sqrt{\hat{p} (1 - \hat{p}) / r}$ where $r$ is the number of replications.  

\subsection{Exponential distribution}
Consider an Exp(1) target distribution, i.e.\ $f(x) = e^{-x}I(x > 0)$.  It is easy to show that $E[X]=1$ and $F^{-1} (q) = \log (1-q)^{-1}$, which we use to evaluate finite sample confidence intervals obtained via fixed-width methods.  We will sample from $f(x)$ using an independence Metropolis sampler with an Exp(1/2) proposal and note this chain is geometrically ergodic \citep{jone:hobe:2001}.  

First, consider estimation of $E[X]$ using each combination of $T_i(\epsilon)$ for $i \in \{ 1, 2, 3 \}$ and $\epsilon \in \{ 0.10, 0.05, 0.02 \}$.  The chain was started from 1 and ran for a minimum of $n^* = 1000$ iterations.  If the stopping criteria was not met, an additional 500 iterations were added to the chain before checking again.  The simulation was repeated for 2000 replications to evaluate the resulting coverage probabilities.

Table~\ref{tab:m1} summarizes the mean and standard deviation of the number of iterations at termination along with the resulting coverage probabilities.  All the coverage probabilities are close to the 0.90 nominal level suggesting all three stopping rules are preforming well.  Note the mean iterations are approximately equal, which is expected since $E[X]=1$ and $\lambda_{\theta} = \Var[X]=1$.

\begin{table}[ht]
\begin{center}
\begin{tabular}{c|cc|cc}
& Length (SD) & $E[X]$ & Length (SD) & $\xi_{.5}$ \\
\hline
$T_1 (0.10)$ & 2.44E3 (4.9E2) & 0.884 & 2.70E3 (5.9E2) & 0.858 \\ 
$T_1 (0.05)$ & 8.89E3 (1.2E3) & 0.894 & 1.01E4 (1.5E3) & 0.881 \\ 
$T_1 (0.02)$ & 5.36E4 (4.7E3) & 0.887 & 6.17E4 (5.4E3) & 0.877 \\ 
\hline
$T_2 (0.10)$ & 2.44E3 (4.8E2) & 0.889 & 5.40E3 (9.4E2) & 0.880 \\ 
$T_2 (0.05)$ & 8.90E3 (1.2E3) & 0.891 & 2.07E4 (2.4E3) & 0.882 \\ 
$T_2 (0.02)$ & 5.35E4 (4.7E3) & 0.887 & 1.29E5 (9.1E3) & 0.883 \\ 
\hline
$T_3 (0.10)$ & 2.45E3 (4.7E2) & 0.888 & 2.79E3 (5.2E2) & 0.865 \\ 
$T_3 (0.05)$ & 8.90E3 (1.2E3) & 0.888 & 1.03E4 (1.3E3) & 0.882 \\ 
$T_3 (0.02)$ & 5.35E4 (4.6E3) & 0.889 & 6.23E4 (5.2E3) & 0.877 \\ 
\end{tabular}
\caption{Summary of coverage probabilities for estimation of $E[X]$ and $\xi_{.5}$ based on 2000 replications and 0.90 nominal level.}
\label{tab:m1}
\end{center}
\end{table}

Next, consider estimation of the median, $\xi_{.5}$, using the same simulation settings.  Table~\ref{tab:m1} summarizes the results from 2000 replications.  Again the results are very close to the 0.90 nominal level, though slightly lower than those for estimating the mean.  Here we have $\xi_{.5} = 0.693$ and $0.5(1-0.5) / e^{\xi_{.5}} = 1$, hence for fixed $\epsilon$ we expect $T_1(\epsilon)$ and $T_3(\epsilon)$ to be similar and $T_2(\epsilon)$ to be larger.

Finally, consider estimating the mean and an 80\% Bayesian credible region simultaneously, which we denote as $\Phi = \left( E[X], \xi_{.1} , \xi_{.9} \right) $.  Due to increased computation time, each chain was run for a minimum of $n^* =10000$ iterations with an additional 5000 added between checks.  The simulation was terminated the first time the length of a confidence interval was sufficiently small for each parameter in $\Phi$.  To adjust for multiplicity, we apply a Bonferonni approach.  Specifically, we set individual confidence intervals to have a coverage probability of $0.90^{1/3} = 0.9655$ resulting in a simultaneous confidence region with coverage probability of at least 0.90.

The simulation was repeated for 2000 replications with each combination of $T_i(\epsilon)$ for $i \in \{ 1, 2, 3 \}$ and $\epsilon \in \{ 0.10, 0.05, 0.02 \}$.  Table~\ref{tab:multi1} summarizes the simulation results.  We can see the individual coverage probabilities improve as $\epsilon$ decreases, especially in the case of $\xi_{.1}$.  For $\epsilon = 0.02$, all the individual coverage probabilities are remarkably close to the nominal level of 0.9655.  Note the observed confidence region coverage probabilities are above the 0.90 nominal level, which is unsurprising due to correlation between parameters in $\Phi$.

\begin{table}[ht]
\begin{center}
\begin{tabular}{c|ccccc}
& Length (SD) & $E[X]$ & $\xi_{.1}$ & $\xi_{.9}$ & Region \\
\hline
$T_1 (0.10)$ & 2.88E4 (3.9E3) & 0.963 & 0.989 & 0.963 & 0.930 \\
$T_1 (0.05)$ & 1.07E5 (9.7E3) & 0.965 & 0.979 & 0.962 & 0.923 \\
$T_1 (0.02)$ & 6.53E5 (3.3E4) & 0.965 & 0.967 & 0.968 & 0.917 \\ 
\hline
$T_2 (0.10)$ & 6.71E4 (5.9E3) & 0.969 & 0.979 & 0.964 & 0.925 \\
$T_2 (0.05)$ & 2.29E5 (1.4E4) & 0.966 & 0.974 & 0.963 & 0.920 \\
$T_2 (0.02)$ & 1.29E6 (5.0E4) & 0.964 & 0.963 & 0.970 & 0.915 \\ 
\hline
$T_3 (0.10)$ & 1.00E4 \hspace{2.5mm} (0) \hspace{2.5mm} & 0.962 & 0.991 & 0.955 & 0.927 \\
$T_3 (0.05)$ & 2.31E4 (2.9E3) & 0.963 & 0.983 & 0.958 & 0.921 \\
$T_3 (0.02)$ & 1.30E5 (9.1E3) & 0.961 & 0.970 & 0.965 & 0.914 \\ 
\end{tabular}
\caption{Summary of coverage probabilities for estimation of $\Phi$ based on 2000 replications.  Individual confidence intervals have a 0.9655 nominal level, resulting in a 0.90 nominal level confidence region.}
\label{tab:multi1}
\end{center}
\end{table}

\subsection{Mixture of bivariate Normals}
Consider a mixture of bivariate Normals ${\bf X} =  \left[ X_1 , X_2 \right]^{T} = p {\bf Y}_1 + (1-p) {\bf Y}_2$, where
\[
{\bf Y}_1= \begin{bmatrix} Y_{11} \\ Y_{12} \end{bmatrix} \sim \N_2 \left( \begin{bmatrix} \mu_{11} \\ \mu_{12} \end{bmatrix}, \begin{bmatrix} \sigma_{11}^2 & 0\\  0& \sigma_{12}^2 \end{bmatrix} \right) \quad \text{ and } \quad
{\bf Y}_2 = \begin{bmatrix} Y_{21} \\ Y_{22} \end{bmatrix} \sim \N_2 \left( \begin{bmatrix} \mu_{21} \\ \mu_{22} \end{bmatrix}, \begin{bmatrix} \sigma_{21}^2 & 0\\  0& \sigma_{22}^2 \end{bmatrix} \right) .
\] 
In this example, we choose $p= 0.25$, $\mu_{11}=1$, $\mu_{12}=10$, $\mu_{21}=2.5$, $\mu_{22}=25$, $\sigma_{11}=0.5$, $\sigma_{12}=5$, $\sigma_{21}=0.7$ and $\sigma_{22}=7$.

We first sample from $f({\bf X})$ with two different component-wise Metropolis random walk algorithms, one with Uniform proposals and another with Normal proposals.  For the Uniform proposals, we apply a $Unif(-3, 3)$ and $Unif(-30, 30)$ random walk for the $X_1$ and $X_2$ dimensions, respectively.   For the Normal proposals, we apply a $N(0, 3^2)$ and $N(0, 30^2)$ random walk for the $X_1$ and $X_2$ dimensions, respectively.  It can be shown that these chains are geometrically ergodic \citep{jarn:hans:2000}.

Consider estimation of $\Phi = \left( E[X], \xi_{.1} , \xi_{.9} \right) $ using fixed-width stopping rules $T_i(\epsilon)$ for $i \in \{ 1, 2, 3 \}$ and $\epsilon \in \{ 0.10, 0.05, 0.02 \}$. We ran the chain for a minimum of $n^* = 5000$ iterations and added 1000 iterations between checking the stopping criteria.  This simulation was repeated for 1000 independent replications.

Table~\ref{tab:ex2metr} summarizes the mean and standard deviation of the number of iterations at termination along empirical coverage probabilities from the Uniform and Normal proposals.  Notice for both samplers, the coverage probabilities improve as $\epsilon$ decreases and are close to the 0.95 nominal level once $\epsilon = 0.02$.  It appears the Metropolis random walk with Normal proposals is mixing faster since the overall simulation effort is substantially lower than that of the Uniform proposals.  This difference in simulation effort illustrates the importance of specifying a good proposal distribution in MCMC simulations.

\begin{table}[ht]
\begin{center}
\begin{tabular}{c|c|cccccc}
\textbf{Uniform} & Length (SD) & $E[X_1]$ & $\xi_{.1, X_1}$ & $\xi_{.9, X_1}$ & $E[X_2]$ & $\xi_{.1, X_2}$ & $\xi_{.9, X_2}$ \\
\hline
$T_1(0.10)$	 & 14,658 (3.4E3) &  0.930 & 0.932 & 0.917 & 0.936 & 0.945 & 0.937\\
$T_1(0.05)$	 & 59,869 (9.1E3) &  0.934 & 0.922 & 0.939 & 0.940 & 0.934 & 0.953\\
$T_1(0.02)$	& 391,566 (3.1E4) & 0.956 & 0.944 & 0.945 & 0.956 & 0.948 & 0.953\\
	\hline
$T_2(0.10)$	& 20,897 (5.0E3) &  0.929 & 0.933 & 0.911 & 0.931 & 0.936 & 0.938\\
$T_2(0.05)$	& 85,401 (1.2E4) &  0.950 & 0.926 & 0.934 & 0.929 & 0.925 & 0.942\\
$T_2(0.02)$	& 556,821(3.9E4) & 0.953 & 0.946 & 0.954 & 0.950 & 0.938 & 0.956\\    
	\hline
$T_3(0.10)$	& 8,827 (1.0E3) &  0.926 & 0.928 & 0.899 & 0.920 & 0.922 & 0.920\\	
$T_3(0.05)$	& 35,733 (2.9E3) &  0.924 & 0.938 & 0.931 & 0.934 & 0.928 & 0.937\\
$T_3(0.02)$	& 233,312 (1.3E4) & 0.954 & 0.955 & 0.959 & 0.948 & 0.958 & 0.956 \vspace*{.5cm} \\
\textbf{Normal} & Length (SD) & $E[X_1]$ & $\xi_{.1, X_1}$ & $\xi_{.9, X_1}$ & $E[X_2]$ & $\xi_{.1, X_2}$ & $\xi_{.9, X_2}$ \\
\hline
$T_1(0.10)$	 & 8,028 (1.5E3) &  0.946 & 0.939 & 0.939 & 0.934 & 0.943 & 0.937\\
$T_1(0.05)$	 & 29,844 (3.7E3) &  0.927 & 0.936 & 0.948 & 0.917 & 0.932 & 0.953\\
$T_1(0.02)$	 & 186,061 (1.3E4) &  0.952 & 0.936 & 0.952 & 0.943 & 0.946 & 0.938\\
	\hline
$T_2(0.10)$	& 11,307 (2.1E3) &  0.949 & 0.933 & 0.940 & 0.940 & 0.944 & 0.943\\
$T_2(0.05)$	& 42,338 (4.6E3) &  0.911 & 0.943 & 0.956 & 0.937 & 0.934 & 0.951\\
$T_2(0.02)$	& 261,741 (1.6E4) &  0.940 & 0.938 & 0.956 & 0.949 & 0.938 & 0.945\\
	\hline
$T_3(0.10)$	& 5,114 (3.2E2) &  0.944 & 0.950 & 0.933 & 0.936 & 0.936 & 0.924\\
$T_3(0.05)$	& 17,654 (1.8E3) &  0.922 & 0.930 & 0.943 & 0.925 & 0.921 & 0.939\\	
$T_3(0.02)$	& 112,626 (7.6E3) &  0.933 & 0.946 & 0.941 & 0.941 & 0.930 & 0.940\\	
\end{tabular}
\caption{Summary of coverage probabilities for estimations of $\Phi$ using a Metropolis random walk with Uniform and Normal proposals based on 1000 replications and a 0.95 nominal level.}
\label{tab:ex2metr}
\end{center}
\end{table}

Next, we consider a Gibbs sampler using the full conditional densities, i.e.\
\begin{align*}
   f_{X_1 | X_2}(x_1 | x_2) & = P_{X_2} Y_{11} + (1-P_{X_2}) Y_{21} \text{ and } \\
   f_{X_2 | X_1}(x_2 | x_1) & = P_{X_1} Y_{12} + (1-P_{X_1}) Y_{22}  \; ,
\end{align*}
where
$$P_{X_2} = \left(1+ {(1-p) \sigma_{12} \over p \sigma_{22}} \exp {\left\{{1\over2} \left( \left({x_2 - \mu_{12} \over \sigma_{12}}\right)^2 - \left({x_2 - \mu_{22} \over \sigma_{22}}\right)^2 \right) \right\}}\right)^{-1},$$
and
$$P_{X_1} = \left(1+ {(1-p) \sigma_{11} \over p \sigma_{21}} \exp {\left\{{1\over2} \left( \left({x_1 - \mu_{11} \over \sigma_{11}}\right)^2 - \left({x_1 - \mu_{21} \over \sigma_{21}}\right)^2 \right) \right\}}\right)^{-1}.$$
Note, $X_1 | X_2 = x_2$ and $X_2 | X_1 = x_1$ are easy to sample from since they are mixtures of Normal random variables.

Table~\ref{tab:ex2gibbs} summarizes the results for the Gibbs sampler.  Notice, the coverage probabilities do not uniformly improve as $\epsilon$ decreases. However, they are all close to the nominal 0.95 level using significantly fewer total iterations, suggesting the Gibbs sampler mixes better than either of the Metropolis random walk samplers.  

As a final comparison, we performed additional simulations via i.i.d.\ sampling (not shown).   The resulting empirical coverage probabilities were similar to using the Gibbs sampler, albeit with slightly fewer iterations.

\begin{table}[ht]
\begin{center}
\begin{tabular}{c|c|cccccc}
\textbf{Gibbs} & Length (SD) & $E[X_1]$ & $\xi_{.1, X_1}$ & $\xi_{.9, X_1}$ & $E[X_2]$ & $\xi_{.1, X_2}$ & $\xi_{.9, X_2}$ \\
\hline
$T_1(0.10)$	 & 1,930 (3.7E2) &  0.941 & 0.940 & 0.937 & 0.954 & 0.958 & 0.927\\
$T_1(0.05)$	 & 5,727 (8.7E2) &  0.946 & 0.958 & 0.941 & 0.942 & 0.945 & 0.940\\
$T_1(0.02)$	 & 31,170 (2.8E3) &  0.935 & 0.945 & 0.961 & 0.937 & 0.937 & 0.944\\
	\hline
$T_2(0.10)$	& 2,465 (5.4E2) &  0.935 & 0.939 & 0.939 & 0.954 & 0.950 & 0.937\\
$T_2(0.05)$	& 7,865 (1.1E3) &  0.950 & 0.959 & 0.943 & 0.955 & 0.954 & 0.952\\
$T_2(0.02)$	& 43,756 (3.6E3) &  0.933 & 0.936 & 0.959 & 0.936 & 0.959 & 0.946\\
	\hline
$T_3(0.10)$	& 1,182 (3.9E2) &  0.929 & 0.936 & 0.942 & 0.936 & 0.936 & 0.924\\
$T_3(0.05)$	& 3,786 (6.2E2) &  0.956 & 0.951 & 0.944 & 0.940 & 0.940 & 0.935\\
$T_3(0.02)$	& 20,289 (2.0E3) &  0.945 & 0.947 & 0.954 & 0.940 & 0.943 & 0.952\\
\end{tabular}
\caption{Summary of coverage probabilities for estimations of  $\Phi$ using a Gibbs sampler based on 1000 replications and a 0.95 nominal level.}
\label{tab:ex2gibbs}
\end{center}
\end{table}


\subsection{Bayesian logistic regression}
Our final example considers the {\it Anguilla} eel data provided in the {\tt dismo} R package \citep[see e.g.][]{elit:leat:2008, hijm:2010}. The data consists of 1,000 observations from a New Zealand survey of site-level presence or absence for the short-finned eel ({\it Anguilla australis}). We selected six out of twelve covariates as in \cite{leat:elit:2008}. Five are continuous variables: SegSumT, DSDist, USNative, DSMaxSlope and DSSlope; one is a categorical variable: Method, with five levels Electric, Spo, Trap, Net and Mixture.

Let $x_i$ be the regression vector of covariates for the $i$th observation of length $k$ and ${\pmb \beta} = \left( \beta_0, \dots, \beta_9 \right)$ be the vector regression coefficients.  For the $i$th observation, suppose $Y_i = 1$ denotes presence and $Y_i = 0$ denotes absence of {\it Anguilla australis}. Then the Bayesian logistic regression model is given by
\begin{align*}
Y_i & \sim Bernoulli(p_i) \; , \\
p_i & \sim {\exp(x_i^{T}{\pmb \beta}) \over 1+\exp(x_i^{T}{\pmb \beta})} \; \text{ and,} \\ 
{\pmb \beta} & \sim N({\pmb 0}, \sigma_{\beta}^2{\bf I}_k) \; ,
\end{align*}
where ${\bf I}_k$ is the $k \times k$ identity matrix. For the analysis, $\sigma_{\beta}^2=100$ was chosen to represent a diffuse prior distribution on ${\pmb \beta}$ \citep{boon:merr:krac:2012}.  Further, we use the {\tt MCMClogit} function in the {\tt MCMCpack} package to sample from the target Markov chain. 

Suppose we are interested in estimating the posterior mean along with an 80\% Bayesian credible interval for each regression coefficient in the model.  Given that we are working with real data, the true values are naturally unknown.  Instead, we ran 1000 independent chains for 1E6 iterations to obtain an accurate estimate, which we treat as the truth (Table~\ref{tab:true3}).

\begin{table}[ht]
\begin{center}
\begin{tabular}{l|ccc}
Variable & $\beta_j$ & $\xi^{(j)}_{.1}$ & $\xi^{(j)}_{.9}$ \\
\hline
Intercept & -10.463 (2.7E-5) &  -12.224 (3.9E-4) & -8.730 (3.7E-4) \\
SegSumT & 0.657 (1.5E-5) & 0.559 (2.1E-5) & 0.757 (2.2E-5) \\
DSDist & -4.02E-3 (3.3E-7) & -6.15E-3 (4.9E-7) & -1.93E-3 (4.4E-7) \\
USNative & -1.170 (7.1E-5) & -1.625 (9.9E-5) & -0.718 (1.0E-4) \\
MethodMixture & -0.468 (6.8E-5) & -0.910 (9.8E-5) & -0.028 (9.8E-5) \\
MethodNet & -1.525 (8.2E-5) & -2.026 (1.2E-4) & -1.035 (1.1E-4) \\
MethodSpo & -1.831 (1.3E-4) & -2.623 (2.2E-4) & -1.798 (1.4E-4) \\
MethodTrap & -2.594 (1.1E-4) & -3.285 (1.8E-4) & -1.937 (1.3E-4) \\
DSMaxSlope & -0.170 (1.1E-5) & -0.244 (1.7E-5) & -0.099 (1.5E-5) \\
USSlope & -0.052 (3.7E-6) & -0.076 (5.5E-6) & -0.028 (5.1E-6) \\
\end{tabular}
\caption{Summary of estimated true values with standard errors for the Bayesian logistic regression example.}
\label{tab:true3}
\end{center}
\end{table}

Consider estimating $\Phi_j = \left( \beta_j, \xi^{(j)}_{.1}, \xi^{(j)}_{.9} \right)$ for $j = 0, \dots, 9$ using fixed-width stopping rules $T_i(\epsilon)$ for $i \in \{ 1, 2, 3 \}$. From the magnitudes in Table~\ref{tab:true3}, it is easy to see a single $\epsilon$ would be problematic for $T_1(\epsilon)$.  Instead, we will specify an $\epsilon$ for each $\Phi_j$ with respect to its magnitude.  Specifically, we choose three simulation settings such that $\pmb{\epsilon}_1$ =  ($1$, $0.01$, $0.001$, $0.1$, $0.1$, $0.1$, $0.1$, $0.1$, $0.01$, $0.01$), $0.5 \pmb{\epsilon}_1$ and $0.2 \pmb{\epsilon}_1$. 

A single $\epsilon$ value for $T_2(\epsilon)$ will also be problematic since there are parameters with very small absolute values (e.g.\ DSDist).  We instead specify an $\epsilon$ for each $\Phi_j$.  In this case, we choose three simulation settings such that $\pmb{\epsilon}_2$ = ($0.1$, $0.1$, $1$, $0.1$, $1$, $0.1$, $0.1$, $0.1$, $0.1$, $1$), $0.5 \pmb{\epsilon}_2$ and $0.2 \pmb{\epsilon}_2$.  

For both $T_1(\epsilon)$ and $T_2(\epsilon)$, it becomes overwhelmingly tedious to specify appropriate $\epsilon$ vectors when the number of parameters becomes large.  However, for the stopping rule $T_3(\epsilon)$ we can use a single $\epsilon$ for the 30 dimensional target parameter vector.  Specifically, we choose three simulation settings such that $\epsilon_3 \in \{0.10, 0.05, 0.02\}$.

For the two larger $\epsilon$ settings, we set $n^* = 10000$ and added 1000 iterations between checks.  For the smallest $\epsilon$ setting, we set $n^* = 1\text{E}5$ and added 10000 iterations between checks due to increased computational demands.  Each simulation setting was repeated 1000 times independently.

Table~\ref{tab:out3} summarizes the empirical coverage probabilities. We can see the coverage probabilities for each stopping rule increase towards the nominal level of 0.95 as $\epsilon$ decreases, suggesting that all the stopping rules perform well.  For high dimensional settings such as this, $T_3(\epsilon)$ provides a distinct practical advantage since a practitioner can specify a single $\epsilon$ value.

\begin{table}[htbp]
\begin{center}
\begin{tabular}{l|ccc|ccc|ccc}
  &\multicolumn{3}{c|}{$T_1(\pmb{\epsilon}_1)$}& \multicolumn{3}{c|}{$T_1(0.5 \pmb{\epsilon}_1)$} & \multicolumn{3}{c}{$T_1(0.2 \pmb{\epsilon}_1)$} \\
\hline
Variable & $\beta_j$ & $\xi^{(j)}_{.1}$ & $\xi^{(j)}_{.9}$ & $\beta_j$ & $\xi^{(j)}_{.1}$ & $\xi^{(j)}_{.9}$ & $\beta_j$ & $\xi^{(j)}_{.1}$ & $\xi^{(j)}_{.9}$ \\
\hline
Intercept & 0.936 & 0.933 & 0.912 &  0.937 & 0.942 & 0.942 &  0.946 & 0.946 & 0.930\\
SegSumT & 0.932 & 0.922 & 0.916& 0.942 & 0.941 & 0.934 & 0.953 & 0.944 & 0.936\\
DSDist & 0.987 & 0.969 & 0.979 & 0.976 & 0.969 & 0.960 & 0.956 & 0.954 & 0.952\\
USNative & 0.927 & 0.929 & 0.917 & 0.939 & 0.933 & 0.943 & 0.948 & 0.939 & 0.944\\
MethodMixture &0.930 & 0.928 & 0.920& 0.946 & 0.948 & 0.938& 0.935 & 0.953 & 0.940\\
MethodNet & 0.946 & 0.922 & 0.936& 0.941 & 0.948 & 0.932 & 0.943 & 0.939 & 0.935\\
MethodSpo & 0.913 & 0.913 & 0.927 & 0.931 & 0.929 & 0.931 & 0.943 & 0.942 & 0.926\\
MethodTrap & 0.928 & 0.906 & 0.937 & 0.938 & 0.930 & 0.927 & 0.941 & 0.947 & 0.947\\
DSMaxSlope & 0.932 & 0.930 & 0.921 & 0.942 & 0.943 & 0.945 &0.953 & 0.958 & 0.951\\
USSlope & 0.921 & 0.928 & 0.935 & 0.951 & 0.927 & 0.954 & 0.957 & 0.952 & 0.962\\
\hline
Length (SD) & \multicolumn{3}{c|}{19,521 (3.8E3)} & \multicolumn{3}{c|}{76,894 (9.5E3)} & \multicolumn{3}{c}{492,910 (3.4E4)} \vspace*{.5cm} \\
  &\multicolumn{3}{c|}{$T_2(\pmb{\epsilon}_2)$}& \multicolumn{3}{c|}{$T_2(0.5 \pmb{\epsilon}_2)$} & \multicolumn{3}{c}{$T_2(0.2 \pmb{\epsilon}_2)$} \\
\hline
Variable & $\beta_j$ & $\xi^{(j)}_{.1}$ & $\xi^{(j)}_{.9}$ & $\beta_j$ & $\xi^{(j)}_{.1}$ & $\xi^{(j)}_{.9}$ & $\beta_j$ & $\xi^{(j)}_{.1}$ & $\xi^{(j)}_{.9}$ \\
\hline
Intercept & 0.928 & 0.938 & 0.915 &  0.950 & 0.948 & 0.947 &  0.945 & 0.949 & 0.938\\
SegSumT & 0.923 & 0.916 & 0.937& 0.953 & 0.955 & 0.948 & 0.944 & 0.947 & 0.947\\
DSDist & 0.985 & 0.968 & 0.975 & 0.970 & 0.958 & 0.958 & 0.956 & 0.955 & 0.947\\
USNative & 0.921 & 0.936 & 0.921 & 0.946 & 0.933 & 0.945 & 0.940 & 0.956 & 0.941\\
MethodMixture &0.941 & 0.938 & 0.933& 0.942 & 0.945 & 0.916& 0.935 & 0.933 & 0.942\\
MethodNet & 0.942 & 0.920 & 0.922& 0.940 & 0.942 & 0.939 & 0.942 & 0.944 & 0.935\\
MethodSpo & 0.919 & 0.901 & 0.924 & 0.936 & 0.923 & 0.937 & 0.947 & 0.956 & 0.947\\
MethodTrap & 0.935 & 0.910 & 0.936 & 0.939 & 0.939 & 0.931 & 0.941 & 0.933 & 0.941\\
DSMaxSlope & 0.937 & 0.942 & 0.916 & 0.948 & 0.942 & 0.950 & 0.942 & 0.954 & 0.955\\
USSlope & 0.935 & 0.933 & 0.930 & 0.949 & 0.936 & 0.941 & 0.949 & 0.944 & 0.943\\
\hline
Length (SD) & \multicolumn{3}{c|}{37,667 (3.5E4)} & \multicolumn{3}{c|}{151,276 (8.9E4)} & \multicolumn{3}{c}{1,161,400 (2.6E5)} \vspace*{.5cm} \\
  &\multicolumn{3}{c|}{$T_3(0.10)$}& \multicolumn{3}{c|}{$T_3(0.05)$} & \multicolumn{3}{c}{$T_3(0.02)$} \\
\hline
Variable & $\beta_j$ & $\xi^{(j)}_{.1}$ & $\xi^{(j)}_{.9}$ & $\beta_j$ & $\xi^{(j)}_{.1}$ & $\xi^{(j)}_{.9}$ & $\beta_j$ & $\xi^{(j)}_{.1}$ & $\xi^{(j)}_{.9}$ \\
\hline
Intercept & 0.932 & 0.944 & 0.929 &  0.943 & 0.950 & 0.943 &  0.943 & 0.954 & 0.934\\
SegSumT & 0.932 & 0.935 & 0.941 & 0.942 & 0.934 & 0.946 & 0.942 & 0.934 & 0.946\\
DSDist & 0.981 & 0.969 & 0.969 & 0.968 & 0.966 & 0.955 & 0.957 & 0.954 & 0.950\\
USNative & 0.939 & 0.942 & 0.923 & 0.941 & 0.948 & 0.954 & 0.942 & 0.943 & 0.940\\
MethodMixture &0.939 & 0.928 & 0.920& 0.947 & 0.943 & 0.933& 0.927 & 0.947 & 0.928\\
MethodNet & 0.929 & 0.922 & 0.931& 0.939 & 0.939 & 0.934 & 0.930 & 0.938 & 0.939\\
MethodSpo & 0.915 & 0.902 & 0.925 & 0.924 & 0.933 & 0.926 & 0.948 & 0.946 & 0.935\\
MethodTrap & 0.930 & 0.909 & 0.920 & 0.941 & 0.937 & 0.933 & 0.939 & 0.935 & 0.948\\
DSMaxSlope & 0.941 & 0.932 & 0.930 & 0.940 & 0.950 & 0.943 & 0.958 & 0.955 & 0.951\\
USSlope & 0.939 & 0.928 & 0.940 & 0.953 & 0.937 & 0.955 & 0.954 & 0.957 & 0.958\\
\hline
Length (SD) & \multicolumn{3}{c|}{24,404 (1.4E3)} & \multicolumn{3}{c|}{78,886 (4.2E3)} & \multicolumn{3}{c}{439,260 (1.7E4)} \\
\end{tabular}
\caption{Summary of coverage probabilities for Bayesian logistic regression example with 1000 independent replicates.  The coverage probabilities have a 0.95 nominal level.}
\label{tab:out3}
\end{center}
\end{table}

To adjust  for multiplicity, we again apply a Bonferonni approach.  We set individual confidence intervals to have a nominal level of $0.80^{1/10}=0.9779$ resulting in simultaneous confidence region with nominal level of at least 0.80.  We only considered estimating the posterior mean of the 10 dimensional vector ${\pmb \beta}$ using $T_3(\epsilon)$ with $\epsilon \in \{ 0.20, 0.10, 0.05, 0.02 \}$. The minimum simulation effort was $n^*=1E5$ iterations with an additional 1000 added between checks. Again, for the smallest $\epsilon$ setting, we set $n^*=1E6$ with an additional 10000 added between checks. The simulation was terminated the first time $T_3(\epsilon)$ was met and repeated 1000 times independently.  

Table~\ref{tab:multiadj} summarizes the simulation results. We can see that, as $\epsilon$ decreases, all the individual coverage probabilities are remarkably close to the nominal level of 0.9779. Note the observed confidence region coverage probabilities approach the nominal level of 0.80 as expected.  However, it is bit surprising how close this is to the nominal 0.80 level given possible correlation among parameters.  To this end, we investigated the correlation between pairs of target parameters.  We found that most pairs have low correlation, except for strong correlation between (Intercept, SegSumT) and moderate correlation between (USNative, USSlope).  Given the lack of correlation, the confidence region coverages are very encouraging.

\begin{table}[ht]
\begin{center}
\begin{tabular}{l|cccc}
 & $T_3(0.20)$ & $T_3(0.10)$ & $T_3(0.05)$ & $T_3(0.02)$ \\
\hline
Variable & $\beta_j$ & $\beta_j$ & $\beta_j$ & $\beta_j$ \\
\hline
	Intercept & 0.959 & 0.975 & 0.976 & 0.973\\
	SegSumT & 0.960 & 0.971 & 0.979 & 0.974\\
	DSDist & 0.995 & 0.989 & 0.993 & 0.979\\
	USNative & 0.948 & 0.978 & 0.970 & 0.973\\
	MethodMixture & 0.950 & 0.973 & 0.967 & 0.968\\
	MethodNet & 0.962 & 0.962 & 0.976 & 0.973\\
	MethodSpo & 0.946 & 0.954 & 0.968 & 0.979\\
	MethodTrap & 0.950 & 0.960 & 0.970 & 0.978\\
	DSMaxSlope & 0.966 & 0.971 & 0.977 & 0.974\\
	USSlope & 0.964 & 0.965 & 0.973 & 0.982\\
\hline
Region & 0.693 & 0.763 & 0.792 & 0.805 \\
\hline
Length (SD) & 10,082(2.7E2) & 29,729(1.8E3) & 100,261(5.2E3) & 583,488(1.9E4)\\
\end{tabular}
\caption{Summary of coverage probabilities for ${\pmb \beta}$ based on $T_3(\epsilon)$ with 1,000 replicates. The coverage probabilities have a 0.9779 nominal level, resulting in a 0.80 nominal level confidence region.}
\label{tab:multiadj}
\end{center}
\end{table}   

\subsection{Discussion}

This paper considers absolute precision, relative magnitude, and relative standard deviation fixed-width stopping rules in the context of MCMC simulations.  Under limited assumptions, we show fixed-width stopping rules obtain a desired coverage probability in an asymptotic sense as $\epsilon$ tends to 0.  Moreover, we illustrate these rules perform well in a variety of finite sample settings provided $\epsilon$ is specified to be small enough.

A practical MCMC stopping rule should be applicable for a large number of parameters since practitioners usually report multiple expectation and quantile estimates.  Unfortunately, choosing a single $\epsilon$ could be problematic for absolute precision and relative magnitude stopping rules.  These stopping rules would be better served by specifying an $\pmb{\epsilon}$ vector, which can be tedious when the number of parameters becomes large.  

Instead, we advocate use of the relative standard deviation stopping rule since it is easy to implement and applicable in multivariate settings without a priori knowledge of the target parameter size.  Simply put, this rule terminates an MCMC simulation when estimates of target parameters are sufficiently accurate relative to their associated posterior standard deviations.  The resulting estimates are approximately $\epsilon^{-1}$ more accurate than their posterior standard deviations.  We recommend using $\epsilon = 0.02$, which provided excellent results in the wide variety of examples considered here.  However, a smaller $\epsilon$ may be appropriate when the accuracy of estimation is critical.

In any MCMC simulation, a key component is choosing a Markov chain that mixes well while sufficiently exploring the state space.  As in the mixture of bivariate Normals, the sampler choice affects the performance significantly in terms of coverage probabilities.  Moreover, the computational effort to achieve a reasonable accuracy varies depending on the sampling scheme.  In practice, the true parameters values are unknown and thus poorly behaved samplers may lead to suspicious inference.  We have offered limited guidance in this direction, but note this is usually the most challenging aspect of an MCMC simulation. An interested reader is directed to \cite{broo:gelm:jone:meng:2010} and the references therein for advice on sampling schemes.

Finally, our examples only consider BM to estimate the asymptotic variance from a CLT since it is the most popular technique and widely available.  Improving the variance estimation step might be possible using alternative methods such as overlapping batch means, spectral variance, or subsampling bootstrap methods \citep{fleg:jone:2010, fleg:2012, fleg:jone:neat:2012}, which are currently available in the {\tt mcmcse} package.



\section*{Acknowledgments}
The authors are grateful to Brian Caffo and Galin Jones for helpful conversations about this paper. 

\begin{appendix}

\section{Proof of Theorem~\ref{thm:T3}} \label{app:T3}
The proof of Theorem~\ref{thm:T3} is very similar to that of Theorem~1 in \cite{glyn:whit:1992}.  A minor modification is necessary due to the relative nature of the stopping rule $T_3 (\epsilon)$.  

Define $z = z_{\delta / 2}$ and $V(n)=2z {\hat{\sigma}_n / \sqrt{n}}+p(n)$, where $p(n) = o(n^{-{1 / 2}})$. Then,
\[
T_3 (\epsilon) = \inf \left\{ n \ge 0 : 2 z \hat{\sigma}_{n} / \sqrt{n} + p(n) \le \epsilon \hat{\lambda}_n \right\} \;
\]
can be denoted as $T_3(\epsilon)= \inf \left\{ n \ge 0 : V(n) \le \epsilon \hat{\lambda}_n \right\}$.  Recall $\sigma^{2}_{\theta} \in (0, \infty)$, then it is easy to verify that 
\begin{equation} \label{eq:p one}
n^{1 / 2}V(n) \rightarrow 2z \sigma_{\theta} > 0 \text{ w.p.1 as } \ n \rightarrow \infty.
\end{equation}
By definition of $T_3 (\epsilon)$, $V(T_3 (\epsilon) - 1) > \epsilon \hat{\lambda}_{T_3 (\epsilon) - 1}$ and there exits a random variable $Z(\epsilon) \in [0, 1]$ such that $V(T_3 (\epsilon) + Z(\epsilon)) \leq \epsilon \hat{\lambda}_{T_3 (\epsilon) + Z(\epsilon)}$.  Further note that $T_3 (\epsilon) \rightarrow \infty$ w.p.1 as $\epsilon \rightarrow 0$ and hence $\hat{\lambda}_{T_3 (\epsilon)} \to \lambda_{\theta}$ w.p.1 as $\epsilon \to 0$.  Then using \eqref{eq:p one} we have
\begin{equation*}
\lim_{\epsilon \rightarrow 0} \sup \epsilon T_3(\epsilon)^{1 / 2} \leq \lim_{\epsilon \rightarrow 0} \sup T_3(\epsilon)^{1 / 2} V(T_3(\epsilon)-1) / \hat{\lambda}_{T_3 (\epsilon) - 1} = 2z \sigma_{\theta} / \lambda_{\theta} \text{ w.p.1.}
\end{equation*}
By a similar argument
\begin{equation*}
\lim_{\epsilon \rightarrow 0} \inf \epsilon T_3(\epsilon)^{1 / 2} \geq \lim_{\epsilon \rightarrow 0} \inf T_3(\epsilon)^{1 / 2} V(T_3(\epsilon)+Z(\epsilon)) / \hat{\lambda}_{T_3 (\epsilon) + Z(\epsilon)} = 2z \sigma_{\theta} / \lambda_{\theta} \text{ w.p.1.}
\end{equation*}
Thus, we have 
\begin{equation}
\label{eq:lim}
\lim_{\epsilon \rightarrow 0} \epsilon T_3(\epsilon)^{1 / 2} = 2z \sigma_{\theta} / \lambda_{\theta} \text{ w.p.1.}
\end{equation}
Given a FCLT at \eqref{eq:clt} holds and $\hat{\sigma}_{n} \to \sigma_{\theta}$ w.p.1 as $n \to \infty$, we have
\begin{equation}
\label{eq:fclt}
\sqrt{n} / \hat{\sigma}_{n} \left( \hat{\theta}_{n} - \theta_{\pi} \right) \stackrel{d}{\to} \; \text{N}(0, 1)
\end{equation}
From \eqref{eq:lim} and \eqref{eq:fclt}, it follows a standard random-time-change argument \citep[p.\ 151 of][]{bill:1995} that 
\begin{equation*}
\sqrt{T_3 (\epsilon)} / \hat{\sigma}_{T_3 (\epsilon)} \left( \hat{\theta}_{T_3 (\epsilon)} - \theta_{\pi} \right) \stackrel{d}{\to} \; \text{N}(0, 1) \text{ as } \epsilon \to 0.
\end{equation*}
Finally, we have
\begin{align*}
\Pr\left( \theta_{\pi} \in C[T_3 (\epsilon)] \right) 
&= \Pr\left(\hat{\theta}_{T_3(\epsilon)} - \theta_{\pi} \in (- z{\hat{\sigma}_{T_3(\epsilon)} / \sqrt{T_3(\epsilon)}}, \ z{\hat{\sigma}_{T_3(\epsilon)} / \sqrt{T_3(\epsilon)}})\right) \\ 
&= \Pr\left(\sqrt{T_3 (\epsilon)} / \hat{\sigma}_{T_3 (\epsilon)} ( \hat{\theta}_{T_3 (\epsilon)} - \theta_{\pi} )) \in (- z, \ z)\right) \to 1 - \delta \text{ as } \epsilon \to 0. 
\end{align*}

\end{appendix}

\setstretch{1}
\bibliographystyle{apalike}
\bibliography{ref}

\end{document}